\documentclass[10pt,twoside,final]{amsart}

\usepackage[english]{babel}
\usepackage{graphicx,epstopdf,epsfig}
\usepackage{amsfonts,epsfig,fancyhdr,graphics,amsmath,amssymb}

\title{Products of involutons  in symplectic groups I: bireflections}

\newtheorem{definition}{Definition}[section]
\newtheorem{theorem}{Theorem}[section]
\newtheorem{proposition}[theorem]{Proposition}
\newtheorem{lemma}[theorem]{Lemma}
\newtheorem{corollary}[theorem]{Corollary}
\newtheorem{remark}[theorem]{Remark}

\newcommand {\ch }{\operatorname{char}}
\newcommand {\SpG }{\operatorname{Sp}}
\newcommand {\PSp }{\operatorname{PSp}}
\newcommand {\GL }{\mathrm{GL}}
\newcommand {\Bahn }{\mathrm{B}}
\newcommand {\rank }{\operatorname{rank}}
\newcommand {\Fix }{\mathrm{Fix}}
\newcommand {\Neg }{\mathrm{Neg}}
\newcommand {\GF }{\mathrm{GF}}
\newcommand {\Idm }{\mathrm{I}}
\newcommand {\Cent }{\mathrm{Cent}}
\newcommand {\rad }{\operatorname{rad}}

\usepackage{ifdraft}
\ifdraft{\usepackage[outer]{showlabels}
	\usepackage{todonotes}}{}
	\ifdraft{\usepackage[outer]{showlabels}
		\usepackage{todonotes} \usepackage{listlbls} \usepackage{datetime}}{}
		
\usepackage[pagebackref,colorlinks,citecolor=red,urlcolor=blue,linkcolor=green, bookmarks=false,hypertexnames=true]{hyperref}
\usepackage{orcidlink}
\usepackage{fancyhdr}

\begin{document}

\bibliographystyle{plain}

\setcounter{page}{1}

\thispagestyle{empty}


\author{Klaus Nielsen}\,\orcidlink{0009-0002-7676-2944}
\email{klaus@nielsen-kiel.de}

\ifdraft{\today \ \currenttime}{\date{November 15, 2024}}
\pagestyle{fancy}
\fancyhf{}
\fancyhead[OC]{Klaus Nielsen}
\fancyhead[EC]{Bireflectionality in symplectic groups}
\fancyhead[OR]{\thepage}
\fancyhead[EL]{\thepage}

\maketitle

\begin{abstract}

We classify bireflectional elements ( products of 2 involutions) in  symplectic groups $\SpG(2n, K)$ over a field $K$. We also classify reversible elements (elements conjugate to their inverses) and bireflectional elements in finite projective symplectic groups $\PSp(2n,q)$.
\end{abstract}


\section{Introduction} \label{intro-sec}

 Let $V$ be a vector space of finite dimension $2n \ge 4$ over a field $K$ of characteristic different from 2. Let $f$ be a nondegenerate symplectic bilinear form on $V$. 

More than 30 years ago, I studied products of involutions in the symplectic group $\SpG(V) = \SpG(V, f)$. I  showed that every element of  $\SpG(V)$ is the product of at most 6 involutions.  Actually, I wanted to show that $\SpG(V)$ is 4-reflectional (a group $G$ is called $k$-reflectional if if  all elements of $G$ are  products of at most $k$ involutions), but I didn't succeed. I only obtained  some partial results:

\begin{theorem} \label{theorem-1}
	\begin{enumerate}
		\item[]
		\item $\SpG(V)$ is 5 reflectional if $K$ is finite.
		\item $\SpG(V)$ is 4 reflectional if $\dim V \equiv 0 \mod 4$  and $(\dim V,|K|) \ne (4,3)$.
	\end{enumerate}
\end{theorem}
For the proof, I needed a classification of the bireflectional elements of $\SpG(V, f)$:

\begin{theorem} \label{theorem-2}
	Let $\varphi \in \SpG(V, f)$. Then   $\varphi$ is bireflectional if and only in a suitable basis of $V$, 
	\[
	\varphi = \left (\begin{array} {cc} A & 0 \\ 0 & A^+ \end{array} \right ),
	f = \left (\begin{array} {cc} 0 & \Idm_n \\ -\Idm_n & 0 \end{array} \right ),
	\]
	where $A^+$ is the transpose inverse of $A$, and $A$ is similar to its inverse.
\end{theorem}
It was already known at that time that $\SpG(V, f)$ is 2-reflectional if $\ch K = 2$.
I also considered the projective symplectic groups:

\begin{theorem} \label{theorem-3}
	\begin{enumerate}
			\item[]
		\item $\PSp(V)$ is 4-reflectional.
		\item $\PSp(V)$ is 3-reflectional if and only if $n$ is odd or $-1 \in K^2 + K^2$.
		\item $\PSp(V)$ is 2-reflectional if and only if $-1 \in K^2$.
	\end{enumerate}
\end{theorem}

I have never published these results so far. But a proof of the 6-reflectionality can be found in the dissertation of H. Röpcke \cite[7.20 Korollar]{Roepcke-1992}. And Bünger \cite[Satz 2.3.5]{Bunger1997} used our theorem \ref{theorem-2} to show that a unitary transformation is bireflectional if and only if it is similar to its inverse and its symplectification is bireflectional.
Meanwhile, there is some new interest in this subject. De la Cruz \cite{Cruz} has shown that $\SpG(V)$ is 4-reflectional if $K$ is the complex field. Ellers and Villa \cite{Ellers_Villa2020} gave a proof of theorem \ref{theorem-1} in the case that $-1 \in K^2$. 

Some years ago, I have also classified the bireflectional elements of the finite projective symplectic groups:

\begin{theorem} \label{theorem-4}
	Let $K = \GF(q)$, where $q \equiv 3 \mod 4$.
	Let $\varphi \in \SpG(V)$. Then   the following are equivalent
	\begin{enumerate}
		\item $\varphi$ is the product of 2 skew-involutions; 
		\item $\varphi$ is reversible;
		\item every elementary divisor $(x\pm 1)^{2t}$ of $\varphi$ occurs with even multiplicity.
	\end{enumerate}
\end{theorem}

We call a linear transformation $\varphi \in \GL(V)$ skew-involutory if $\varphi^2 = -1$.

\begin{theorem} \label{theorem-5}
	Let $K = \GF(q)$, where $q \equiv 3 \mod 4$.
	Let $\varphi \in \SpG(V)$. Then   the following are equivalent
	\begin{enumerate}
		\item $\varphi$ is the product of an involution and a skew-involution;
		\item $\varphi$ is conjugate to $-\varphi^{-1}$,
		and every elementary divisor $(x^2 + 1)^{2t}$ of $\varphi$ occurs with even multiplicity.
	\end{enumerate}
\end{theorem}

It follows that $\PSp(2n, q)$ contains reversible elements that are not bireflectional if and only if $q \equiv 3 \mod 4$ and $2n \ge 8$:

\begin{corollary} \label{cor}
	Let $K = \GF(q)$, where $q \equiv 3 \mod 4$.
	Let $\varphi \in \SpG(V)$. Then  the image of $\varphi$ in $\PSp(V)$ is reversible but not bireflectional if and only if
	\begin{enumerate}
		\item $\varphi$ is conjugate to $-\varphi^{-1}$,
		\item $\varphi$ has an elementary divisor $(x^2 + 1)^{2t}$ of odd multiplicity, and
		\item $\varphi$ has an elementary divisor $(x - 1)^{2t}$ of odd multiplicity.
	\end{enumerate}
\end{corollary}

Recently, C. de Seguins Pazzis proved our theorems \ref{theorem-2} (see  \cite[Theorem 1.1]{Pazzis-2024a})  and  \ref{theorem-1} (see \cite[Theorem 1.5]{Pazzis-2024b}  and \cite[Theorem 1.1]{Pazzis-2025}). Moreover, he showed that all symplectic groups are 5-reflectional.

So I think it is time to present the old short proof of theorem \ref{theorem-2}. 
We also prove theorems \ref{theorem-4} and \ref{theorem-5}.
In a following paper we deal with theorems \ref{theorem-1} and \ref{theorem-3} and prove that $\SpG(V)$ is 5-reflectional.

\section{Preliminaries}
 
 \subsection{Notations}
 \mbox{}

 For a linear mapping $\varphi$ of $V$ let $\Bahn^j(\varphi)$ denote the image and $\Fix^j(\varphi)$ the kernel of $(\varphi -1)^j$. Put $\Fix^{\infty}(\varphi) = \bigcup_{j\ge 1} \Fix^j(\varphi)$ and $\Bahn^{\infty}(\varphi) = \bigcap_{j\ge 1}
 \Bahn^j(\varphi)$. The space $\Bahn(\varphi) := \Bahn^1(\varphi)$ is the path or residual space of $\varphi$.
 And $\Fix(\varphi) := \Fix^1(\varphi)$ is the fix space of $\varphi$. Further, let $\Neg^j(\varphi) =  \Fix^j(-\varphi)$. The space $\Neg(\varphi) = \Neg^1(\varphi)$ is the negative space of  $\varphi$.
 We call $\varphi$ a big transvection if $\Bahn(\varphi) \le \Fix(\varphi)$. Hence a big transvection  is unipotent of index 2, and a transvection is a big transvection with path dimension one.
 
 By $\mu(\varphi)$ we denote the minimal polynomial of 
 $\varphi$.

 \subsection{The symplectic matrix group}
 \mbox{}
 
 Usually, we work basis-free, but sometimes it is more convenient to compute with matrices. 
A matrix $P$ is an automorph of a matrix $G$ if $PGP' = G$, where $P'$ is the transpose of $P$.
The skew-sum $[A \setminus B]$ of 2 matrices $A$ and $B$ is the matrix 
\[
\left (\begin{array} {cc} 0 & B\\ A & 0 \end{array} \right ).
\]
By $\SpG(2n, K)$ we denote the standard symplectic group, the set of all automorphs of the matrix $G := [-\Idm_n \setminus \Idm_n]$.

 \subsection{Orthogonally indecomposable transformations and Wall form}
 \mbox{}

According to Huppert \cite[1.7 Satz]{Huppert-1980a}, an orthogonally indecomposable transformation of $\SpG(V)$ is either\footnote{Our type enumeration follows Huppert \cite{Huppert-1980a}. In \cite{Huppert-1990} Huppert uses  a different  enumeration.}
\begin{enumerate}
	\item bicyclic with elementary divisors $e_1 = e_2 =(x \pm 1)^m$
	(type 1) or
	\item indecomposable as an element of $\GL(V)$ (type 2) or
	\item cyclic with minimal polynomial $(qq^*)^t$, where $q$ is irreducible and prime to its reciprocal $q^*$, where $q^*(x) = q(0)^{-1}x^{\partial q}q(x^{-1})$ ( type 3).
\end{enumerate}

\begin{remark}
	Let $\varphi \in \SpG(V)$. Then 
	$Vf(\varphi)^{\bot} = \ker f^*(\varphi)$for all polynomials $f(x) \in K[x]$. In particular, $\Bahn^j(\varphi)^{\bot} = \Fix^j(\varphi)$ and $V = \Bahn^{\infty}(\varphi) \perp \Fix^{\infty}(\varphi)$.	
\end{remark}

\begin{definition}
	Let $\varphi \in \SpG(V)$. We say that $\varphi$ is
	hyperbolic if $V$ is the direct sum of 2 totally isotropic $\varphi$-invariant subspaces.
\end{definition}

In \cite[2.3 Satz]{Huppert-1980a} and \cite[2.4 Satz]{Huppert-1980a}, Satz], Huppert proved the following 

\begin{proposition} \label{prop-1}
	Let $\varphi \in \SpG(V)$ be orthogonally indecomposable  of type 1. Then $\dim V \equiv 2 \mod 4$, and $\varphi$ is hyperbolic.
\end{proposition}

\begin{lemma}            \label{WALL1}
	Let $K$ be finite. Let $\varphi \in \SpG(V)$ be bi- and homocyclic.
	If $\varphi^2$ is fixfree, then $\varphi$ is hyperbolic.
\end{lemma}
	
\begin{proof}
	This follows immediately from a result of Wall; see \cite[2.6 example]{Wall}	on p. 36.
\end{proof}
		
\begin{definition}[Wall form]
	Let $\varphi \in \SpG(V)$. Then 
	$\omega_{\varphi}(u(1- \varphi), w (1-\varphi)) = f(u, w (1-\varphi))$
	defines a nondegenerate bilinear form, the Wall form of $\varphi$, on the path $\Bahn(\varphi)$ of $\varphi$. Let $\Theta(\varphi)$ denote the discriminant of $\omega_{\varphi}$.
\end{definition}

\begin{remark}            \label{WALL2}
	Let $\varphi \in \SpG(V,f)$. Then 
	$\omega_{\varphi}(u,w) - \omega_{\varphi}(w,u) = f(u, w)$.
\end{remark}

\begin{remark}            \label{WALL3}
Let $\varphi \in \SpG(V,f)$  be a big transvection. Then $\omega_{\varphi}$ is symmetric, and $\varphi$ is conjugate to matrix 
\[
P = \left (\begin{array} {cc} \Idm_m  & S  \\ 0 &  \Idm_m\end{array} \right ) \in \SpG(2n, K),
\]
where $S = 0 \oplus T$ and $T$ is congruent to $\omega_{\varphi}$.
\end{remark}

It is  known that 2 cyclic unipotent isometries are conjugate if and only if their  Wall forms have the same discriminant. For the convenience of the reader, we provide a short proof.

\begin{lemma}            \label{WALL4}
	 Let $\varphi \in \SpG(V,f)$ be cyclic with minimal polynomial $(x-1)^{2n}$. Let $\theta \in \Theta(\varphi)$.
	Then $\omega_{\varphi}$ is congruent to a unique  upper antitriangular matrix  $A=(a_{i,j})$ satisfying 
	\begin{enumerate}
		\item $a_{i,i} = 0$ if $i \le  n-2$,
		\item $a_{i,j} = -a_{j+1,i-1}$,
		\item $a_{i,j} - a_{j,i} = a_{i+1,j}$,
		\item $a_{i,j} = 0$ if $i +j \ge  2n+1$,
		\item $a_{j, 2n-j} = (-1)^{j-1} a_{1, 2n-1}$
		\item $a_{n,n} = (-1)^{n+1} \theta$.
	\end{enumerate}
\end{lemma}

\begin{proof}
	There exists a vector $u \in V - \Bahn(\varphi)$ such that $f(u (1-\varphi)^j ,u(1-\varphi)^{j+1}) = 0$ for $0 \le j \le n-2$:
	This is obviously true for $n = 1$. So let $n \ge 2$.
	By induction, there exist a vector $w \in \Bahn(\varphi)- \Bahn^2(\varphi)$ such that $f(w (1-\varphi)^j,w (1-\varphi)^{j+1}) = 0$ for $0 \le j \le n-3$.
	Let $w = v - v \varphi$ for some vector $v \in V$. There exists a vector
	$z \in \Fix^2(\varphi) - u(\varphi - \varphi^{-1})^{\perp}$. Then 
	$f(v + \lambda z, v \varphi + z \varphi) = f(v, v \varphi) + f(\lambda z, v(\varphi - \varphi^{-1}) = 0$ for some $\lambda \in K$.
	Put $u = v + \lambda z$. Then $f(u (1-\varphi)^j, u (1-\varphi)^{j+1} =
	f(w (1-\varphi)^{j-1}, w (1-\varphi)^j = 0$ for $j \ge 2$ and
	 $f(u (1-\varphi), u (1-\varphi)^2) = f(w + \lambda z(1-\varphi), w (1-\varphi)) = 0$.
	
	Now we compute
	$\omega_{\varphi}$ in the basis $ u_1 = u(1-\varphi), \dots, u_{2n-1} = u(1-\varphi)^{2n-1}$ of $\Bahn(\varphi)$.
	For $i, j \ge 1$ define $a_{i,j} = \omega_{\varphi}(u_i, u_j) = f(u_{i-1},u_j)$.
	Then 
	\begin{enumerate}
		\item $a_{i,i} = 0$ if  $i \le  n-2$.
		\item $a_{i,j} = f(u_{i-1},u_j) = -f(u_j, u_{i-1}) = -a_{j+1,i-1}$. 
		\item $a_{i,j} - a_{j,i}  = \omega_{\varphi}(u_i, u_j) - \omega_{\varphi}(u_j, u_i) = f(u_i, u_j)= a_{i+1,j}$ by \ref{WALL2}.
		\item $a_{i,j} = 0$ if $i+j \ge 2n+1$, as $u_{i-1}(1-\varphi^{-1})^j = 0$ if $i+j \ge 2n+1$.
		\item $a_{j, 2n-j} = f(u_{j-1}, u_{2n-j}) = f(u, u_{n-j}(1-\varphi^{-1})^{j-1})
		= (-1)^{j-1} f(u, u_{2n-1}) = (-1)^{j-1} a_{1, 2n-1}$, as $u_{2n-j}(1-\varphi^{-1})^{j-1} \in \Fix(\varphi)$.
		\item $\Theta(\varphi) = (-1)^{n+1} a_{n,n}$.
	\end{enumerate}
	
	We show by induction that if $B =(b_{i,j})$ is  another matrix satisfying (1) - (6), then $B=A$:
	Clearly, $B=A$ if $n = 1$.
	So we may assume that $a_{i,j} = b_{i,j}$ for 
	$2 \le i,j \le 2n-2$. Then
	
	\begin{enumerate}
		\item[(i)]  $a_{1,1} = 0 = b_{1,1}$ by (1),
		\item[(ii)]  $a_{2,1} = 0 = b_{2,1}$ by (2),
		\item[(iii)]  $a_{1,2} = a_{2,2} + a_{2,1} =  b_{2,2} + b_{2,1} = b_{2,1}$ by (ii) and (3)
		\item[(iv)]  $a_{i,1} = - a_{2,i-1} = - b_{2,i-1} = b_{i,1} $ for $3 \le i \le 2n-1$ by (2),
		\item[(v)]  $a_{1,j} = a_{2,j} + a_{j,1} = b_{2,j} + b_{j,1} = b_{1,j} $ for $3 \le j \le 2n-2$ by (3),
		\item[(vi)]  $a_{1,2n-1} = a_{2n-1, 1} = b_{2n-1, 1} = b_{1,2n-1}$ by (5) and (iv),
		\item[(vii)]  $a_{2n-1, j} = 0 = b_{2n-1, j}$ and $a_{j, 2n-1} = 0 = b_{j, 2n-1}$ for $j \ge 2$ by (4).
	\end{enumerate}
\end{proof}

\begin{corollary} \label{WALL5}
	 Let $\varphi, \psi  \in \SpG(V,f)$ be unipotent and cyclic. Then $\varphi$ and $\psi$ are conjugate iff
	$\Theta(\varphi) = \Theta(\psi)$.
\end{corollary}

\begin{proof}
	By a result of Wall \cite[Theorem 1.3.1]{Wall}, $\varphi$ and $\psi$ are conjugate if and only if $\omega_{\varphi}$ and $\omega_{\psi}$ are congruent.
	Apply \ref{WALL4}.
\end{proof}

\begin{corollary}            \label{WALL6}
	Let $\varphi \in \SpG(V,f)$ be unipotent and cyclic. Let $\widehat{\varphi}$ the symplectic transformation induced by $\varphi$ on $\widehat{V} := \Bahn(\varphi)/\Fix(\varphi)$. Then  
	$\Theta(\varphi) = -\Theta(\widehat{\varphi})$.
\end{corollary}

\section{Products of 2 symplectic involutions: Proof of theorem \ref{theorem-2}}
\mbox{}

\begin{lemma}                                                \label{SYMPINV0}
	Let $\varphi  \in \SpG(V)$ be the product of two 
	involutions   $\sigma, \tau  \in \SpG(V)$. Then
	$T\sigma, T\tau \subseteq
	T^{\perp}$ for every  $\varphi$-cyclic subspace $T$.
\end{lemma}

\begin{proof}
	Let $v \in V$. Then
	\[
	f(v\varphi^i \sigma, v\varphi^{i+2k})
	=  f(v\varphi^i \sigma \varphi^{-k}, v\varphi^{i+k})
	=  f(v\varphi^{i+k}\sigma, v\varphi^{i+k}) = 0,
	\]
	as $\sigma$ is alternating. Since $\sigma \varphi = \tau$
	is alternating, we also have
	\[
	f(v\varphi^i \sigma, v\varphi^{i+2k-1}) =
	f(v\varphi^i \sigma \varphi^{-k+1}, v\varphi^{i+k}) =
	f(v\varphi^{i+k}\sigma \varphi, v\varphi^{i+k}) = 0.
	\]
\end{proof}

We prove a symplectic analogon of a result in orthogonal groups \cite[Proposition]{KN-1987a}:

\begin{lemma}                                                \label{SYMPINV1}
	Let $\varphi \in \SpG(V)$, and let $\sigma \in \SpG(V)$ be
	an involution inverting $\varphi$.  There
	exists a decomposition
	$V = V_1 \perp \dots \perp V_s$ of $V$ into an orthogonal sum
	of $\langle \varphi, \sigma \rangle$-invariant subspaces such that
	$V_i = U_i \oplus W_i$, where $U_i$ and $W_i$ are
	totally degenerate, $\varphi$-cyclic  and $\sigma$-invariant.
\end{lemma}

\begin{proof}
	Let $\mu_{\varphi}$ be the minimal polynomial of $\varphi$. If
	$\mu_{\varphi} = pr$, where $p,r \in K[x]$ are selfreciprocal
	polynomials with $\gcd(p,r) = 1$, then
	$V = \ker p(\varphi) \perp \ker r(\varphi)$, and $\ker p(\varphi)$ and
	$\ker r(\varphi)$ are $\sigma$-invariant.
	So by induction, we
	may assume that either $\mu_{\varphi} = p^t$ is the power of an
	irreducible selfreciprocal polynomial $p$ or
	$\mu_{\varphi} = (pp^*)^t$, where $p$ is irreducible and prime to
	$p^*$. Let $U$ be the $\varphi$-cyclic subspace  generated by some
	$u \in [\Neg(\sigma) \cup \Fix(\sigma)] - [V p(\varphi)
	\cup V p(\varphi^*)]$.  Let $W$ be the  $\varphi$-cyclic subspace generated
	by some $w \in [\Neg(\sigma) \cup
	\Fix(\sigma)] - [(u p(\varphi)^{t-1})^{\perp} \cup
	(u p^* (\varphi)^{t-1})^{\perp}]$. Then $U$ and $W$ are $\sigma$-invariant
	and $U \oplus W$ is regular.
	By \ref{SYMPINV0}, $U$ and $W$ are totally
	degenerate.
\end{proof}


\section{Products of 2 skew-involutions}
\mbox{}

\begin{lemma}                                                  \label{SKEWINV0}
	Let  $\lambda \in K$. Let  $\varphi \in \GL(V)$ be cyclic and similar to $\lambda \varphi^{-1}$. Assume further that $\varphi^2$ has no elementary divisor $(x-\lambda)^t$ of odd degree $t$. Then $\dim V = 2m$ is even, and $\varphi$ is similar to a matrix
	\[
	P = \left (\begin{array} {cc} 0 & \Idm_m  \\ -\lambda \Idm_m & D \end{array} \right ).
	\]
\end{lemma}

\begin{proof}
	We may assume  that $K$ is algebraicly closed.
	Let $\mu^2 = \lambda^{-1}$, and put $\psi = \mu \varphi$. Then $\psi$ is similar to its inverse.  Hence $\psi$ is the product of 2 
	involutions $\sigma$ and $\tau$. This is a result usually attributed to Wonenburger \cite[Theorem 1]{Wonenburger-1966} and Dokovi\'c \cite[Theorem 1]{Dokovic-1967}, but it was already known at the beginning of the 20th century for complex matrices; see e.g.
	Jackson \cite[p. 480]{Jackson1909}. Frobenius \cite[§3]{Frobenius-1910} observed that if $P^A = P^{-1}$ and $A^2 = f(A^2)^2$ for some polynomial $f$, then $A f(A^2)^{-1}$ is an involution inverting $P$. And in 1896  he 
	 had already shown that a nonsingular complex matrix $M$ has a square root that is a polynomial in $M$\footnote{Über die cogredienten Transformationen der bilinearen Formen. Berlin Sitzb. (1896), p. 4}. This results holds true for all algebraicly closed fields of characteristic different from two. To prove this we may assume that $M$ is homocyclic: If 
	 \[
	 M = \left (\begin{array} {cc} A & B  \\ 0 & D \end{array} \right )
	 \]
	and $M = f(M)^2$, then $A = f(A)^2$ and $D = f(D)^2$. 
	Moreover, we may assume that $M$ is cyclic. But then $\Cent(M) = K[M]$. This is another result  of Frobenius\footnote{stated without proof in 'Über lineare Substitutinen und bilineare Formen', Crelle 84 (1878) } \cite[(8.)]{Frobenius-1910}.

	The spaces $\Neg^{\infty}(\psi)$ and $\Bahn^{\infty}(-\psi)$ are $\sigma$-invariant and $\tau$-invariant.
	So we may assume that $\Neg(\psi) = 0$.
	
	Then  $\Fix(\sigma) \cap \Neg(\tau) = 0$ and $\Fix(\tau) \cap \Neg(\sigma) = 0$. Further, we may assume that $\Fix(\sigma) \cap \Fix(\tau) = 0$ as $\psi$ is cyclic. Then 
	$\Fix(\sigma) \cap \Fix(\sigma) \tau = 0$. Hence
	$V = \Fix(\sigma) \oplus \Fix(\sigma) \tau$, and in a suitable basis,
	\[
	\sigma = \left (\begin{array} {cc} \Idm_m & 0 \\ C & -\Idm_m \end{array} \right ),
	\tau = \left (\begin{array} {cc} 0 & \Idm_m \\ \Idm_m & 0 \end{array} \right ), 
	\psi = \left (\begin{array} {cc} 0 & \Idm_m  \\ -\Idm_m & C \end{array} \right ).
	\]
	Hence $\varphi$ is similar to
	\[
	 \left (\begin{array} {cc} 0 & \mu^{-1} \Idm_m  \\ - \mu^{-1}\Idm_m & \mu^{-1} C \end{array} \right )
	 \sim \left (\begin{array} {cc} 0 & \Idm_m  \\ -\lambda \Idm_m & \mu^{-1} C \end{array} \right ) =: Q.
	\]
	Now $Q \oplus \lambda Q^{-1} =  \mu^{-1} C \oplus \mu^{-1} C$. Hence $\mu^{-1} C$ is similar to a matrix $D$ in $K^{n,n}$.
\end{proof}

\begin{remark}  
	The Dickson transform of a monic polynomial $f$ of degree $n$ is the polynomial $f^{\mathcal{D}}(x) = f(x+x^{-1})x^n$. Let
	\[
	P =
	\left (\begin{array} {cc} 0 & \Idm_n\\ -\Idm_n & D \end{array} \right ).
	\]
	Then the invariant factors of $P$ are the Dickson transforms
	of the invariant factors of $D$. 
\end{remark}

\begin{lemma}                              \label{SKEWINV1} 
	Let $K$ be finite.
	Let $C \in \GL(m,K)$ be irreducible. Then $K[C]$ is a finite field.
	Hence there exists matrices
	$A, B \in K[C]$ such that $A^2 + B^2 + ABC = -\Idm_m$. Put $D = B+AC$.
	Then
	\[
	\left (\begin{array} {cc} 0 & \Idm_m \\ -\Idm_m & C \end{array} \right ) =
	\left (\begin{array} {cc} -D & A \\ A - CD & D \end{array} \right )
	\left (\begin{array} {cc} A & B \\ D & -A \end{array} \right )
	\]
	is a product of two skew-involutions.
\end{lemma}

\begin{lemma}                                                 \label{SKEWINV2}
	Let $K$ be finite.
	Let $\varphi \in \GL(V)$ be cyclic with minimal polynomial $p^t$, where
	$p$ is irreducible of even degree. Let $\varphi$ be reversible. Then
	$\varphi$ is a product of 2 skew-involutions.
\end{lemma}

\begin{proof}
Let $V = \langle u \rangle_{\varphi}$, and let $\varphi_\mathrm{S}$ the semisimple factor in the Jordan-Chevalley decomposition of $\varphi$.
Put $W = \langle u \rangle{\varphi_\mathrm{S}}$.
		
By \ref{SKEWINV0} and \ref{SKEWINV1}, the restriction of
$\varphi_\mathrm{S}$ on $W$ is inverted by a skew-involution.
Hence  there exists a polynomial $g \in K[x]$ such that $u g(\varphi_\mathrm{S})  g(\varphi_\mathrm{S}^{-1}) = -u$.
		
There exist a polynomial $h$ such that $g(\varphi_\mathrm{S}) =  h(\varphi)$.
Let  $\alpha \in \GL(V)$ with $\varphi^{-1} = \varphi^{\alpha}$. 
Then 
$\varphi_\mathrm{S}^{-1} = \varphi_\mathrm{S}^{\alpha}$ and 
$h(\varphi^{-1}) = h(\varphi^{\alpha}) = h(\varphi)^{\alpha} = g(\varphi_\mathrm{S})^{\alpha} = g(\varphi_\mathrm{S}^{\alpha}) = g(\varphi_\mathrm{S}^{-1})$.

Put $w = u h(\varphi)$.
Define the linear transformation $\eta$ by
$u \varphi^j \eta = w \varphi^{-j}$. Then for $f \in K[x]$, $u f(\varphi)\eta = w f(\varphi^{-1})$ and
$u \varphi^j \eta^2 =  w \varphi^{-j} \eta = u h(\varphi) \varphi^{-j}\eta = w h(\varphi^{-1}) \varphi^j = u h(\varphi) h(\varphi^{-1}) \varphi^j = -u \varphi^j$.
		
Further $u \varphi^j (\varphi \eta)^2 = w \varphi^{-j-1} \varphi \eta = w \varphi^{-j} \eta = -u\varphi^j$. Hence $\eta$ and $\varphi \eta$ are skew-involutions.
\end{proof}

\begin{lemma}                                             \label{SKEWINV3}
	Let $\varphi \in \SpG(V, f)$ be cyclic.
	Let  $\sigma \in \GL(V)$ be an involution or a  skew-involution inverting $\varphi$. 
	\begin{enumerate}
		\item If $\sigma $ is an involution, then $\sigma$ is skew-symplectic.
		\item If $\sigma $ is a skew-involution, then $\sigma$ is symplectic.
	\end{enumerate}
\end{lemma}

\begin{proof}
	Let $V = \langle u \rangle_{\varphi}$, and put $w = u \sigma$. 
	Then
	$f(u \varphi^i \sigma, w \varphi^j \sigma) = f(u \sigma \varphi^{-i}, w \sigma \varphi^{-j}) = f(u \sigma \varphi^j, w \sigma\varphi^i) = - f(w \sigma \varphi^i, u \sigma \varphi^j) = -\epsilon  f(u \varphi^i, \varphi^j)$,
	 where $\sigma^2 = \epsilon$.
\end{proof}

\begin{lemma}                                             \label{SKEWINV4}
	Let $P \in \SpG(2n, K)$ be hyperbolic.
	Then $P$ is the product of 2 symplectic skew-involutions.
\end{lemma}
	
\begin{proof}
		We may assume that 
		\[
		P =
		\left (\begin{array} {cc} Q & 0\\ 0 & Q^+ \end{array} \right ).
		\]
		Now $Q = ST$ is the product of 2 symmetric matrics $S$ and $T$. This is a well-known fact; see e.g. \cite[Satz 5a]{Shoda1929}. (In fact, all matrices conjugating a cyclic matrix into its transpose must already be symmetric. To see this, one may assume that $K$ is algebraicly closed and use Frobenius' work). Now 
		\[
		P =
		\left (\begin{array} {cc} 0 & S\\ -S^{-1} & 0 \end{array} \right )
		\left (\begin{array} {cc} 0 & -T^{-1}\\ T & 0 \end{array} \right )
		\]
		is the product of 2 symplectic skew-involutions.
\end{proof}

\begin{remark}                                    
	Let $\varphi \in \SpG(V)$ be unipotent and cyclic. If $\varphi$ is conjugate to its inverse, then $-1 \in K^2$. 
\end{remark}

\begin{proof}
	Let $\alpha \in \SpG(V)$ and $\varphi^{\alpha} = \varphi^{-1}$. Let $\sigma$ be an involution reversing $\varphi$.  Then $\alpha \sigma$ centralizes $\varphi$. 
	Hence $\alpha \sigma$ is a polynomial in $\varphi$, 
	and $\alpha \sigma$ has a single eigenvalue $\lambda$.
	By \ref{SKEWINV3}, 
	$\sigma$ must be skew-symplectic.
	Hence  $\alpha \sigma$ is skew-symplectic. It follows that $\lambda = - \lambda^{-1}$.
\end{proof}

\begin{lemma}                                     \label{SKEWINV5}
	Let $P \in \SpG(2n, K)$
	be conjugate to its inverse. If $-1 \not \in K^2$, then
	 $\dim \Bahn^t(P)$ is even.
\end{lemma}

\begin{proof}
	By induction, we may assume that $t=1$. 
	Let 
	\[
	Q =\left (\begin{array} {cc} A & B\\ C & D \end{array} \right),
	P = \left (\begin{array} {cc} \Idm_n & 0\\ S & \Idm_n \end{array} \right)
	\]
	with $P^Q = P^{-1}$. From 
	\[
	\left (\begin{array} {cc} A & B\\ C & D \end{array} \right)
	\left (\begin{array} {cc} \Idm_n & 0\\ S & \Idm_n \end{array} \right) =
	\left (\begin{array} {cc} A+BS & B\\ C+DS & D \end{array} \right),
	\]
	\[
	\left (\begin{array} {cc} \Idm_n & 0\\ -S & \Idm_n \end{array} \right)
	\left (\begin{array} {cc} A & B\\ C & D \end{array} \right)=
	\left (\begin{array} {cc} A & B\\ C-SA & D-SB \end{array} \right)
	\]
	it follows that $SA=-DS, SB=0=BS$.  Since $AD' - BC' = \Idm_n$, we obtain  $S = -DSD'$. Further, $S$ is symmetric, hence 
	$\dim \Bahn(P) = \rank S$ must be even.
\end{proof}

\begin{lemma}                                     \label{SKEWINV6}
	Let $\varphi \in \SpG(V,f)$ be bicyclic with elementary divisors $(x-1)^{2m}$.
	Then $g_{\varphi}(a,b) = f(a(1-\varphi)^{2m-1},b)$
	defines a symmetric bilinear form on $V$ with $\rad g_{\varphi} = \Bahn(\varphi)$. Let $\widehat{V} = \Bahn(\varphi)/\Fix(\varphi)$ and $\widehat{\varphi} = \varphi|_{\widehat{V}}$.
	\begin{enumerate}
	  \item If $\varphi$ is hyperbolic, then $V/\rad g_{\varphi}$ is hyperbolic.
      \item $V/\rad g_{\varphi}$ is hyperbolic if and only if $\widehat{V}/\rad g_{\widehat{\varphi}}$ is hyperbolic.
		\end{enumerate}
\end{lemma}

\begin{proof}
	Clearly, $V/\rad g_{\varphi}$ is hyperbolic if $\varphi$ is hyperbolic.
	We prove 2: 
	
	Let $a \in V$ and $\widehat{a} = a (\varphi^2 - 1)$.  Then
	
	$f(\widehat{a} (\varphi^2 - 1)^{2m-3}, \widehat{a}) = f(a (\varphi^2 - 1)^{2m-2}, a (\varphi^2 - 1)) = f(a (\varphi^2 - 1)^{2m-2} (\varphi^{-2} - 1), a ) = -f(a(\varphi^2 -1)^{2m-1}, a)$. Hence $a$ is isotropic if and only if $\widehat{a} + \Fix(\varphi)$ is isotropic.
\end{proof}

\begin{lemma}                                                                            \label{SKEWINV7}
	Let  $P \in \SpG(4, K)$ be cyclic with minimal polynomial $(x^2+1)^2$.
	Then $-P^2$ is a big transvection with Wall form $\Idm_2$.
\end{lemma}

\begin{proof}
	Let 
	\[
	H =\left (\begin{array} {cc} 0 & 1\\ -1 & 0 \end{array} \right),
	P = \left (\begin{array} {cc} H & B\\ 0 & H \end{array} \right),
	T = \left (\begin{array} {cc} -\Idm_2 & HB + BH \\ 0 & -\Idm_2 \end{array} \right).
	\]
	Then $B$ has trace zero, $HB+BH = \lambda \Idm_2$ for some $\lambda \ne 0$, and $P^2 = T$.
	Apply \ref{WALL3}.
\end{proof}

Using \ref{SKEWINV6} and \ref{SKEWINV7}, we obtain
\begin{corollary}                                                                            \label{SKEWINV8}
	Let $K = \GF(q)$ be a finite field with $q \equiv 3 \mod 4$, 
	and let $\varphi \in \SpG(V)$ be cyclic with minimal polynomial $(x^2+1)^n$.
	Then $\varphi^2$ is hyperbolic if and only if $n$ is odd.
\end{corollary}

\begin{proof}[Proof of theorem \ref{theorem-4}]  
	$1 \rightarrow 2$: Apply  \ref{SKEWINV5}. 
	The implication $3 \rightarrow 1$ is  clear.
	
	$2 \rightarrow 3$: 
	We may assume that $\varphi$ is either orthogonally indecomposable
	or bicyclic with elementary divisors $(x-1)^{2t}$.
	
	If $\varphi$ is  orthogonally indecomposable of type $2$ with minimal polynomial $p(x)^t$, where $p$ is irreducible of degree $\ge 2$, then by \ref{SKEWINV2}, $\varphi$ is reversed by a skew-involution $\eta \in \GL(V)$, and by \ref{SKEWINV3}, $\eta \in \SpG(V)$.
	
	If $\varphi$ is orthogonally indecomposable of type 1 or 3, then $\varphi$ is hyperbolic. And if $\varphi$ is hyperbolic we can apply \ref{SKEWINV4}.
	
	It remains to consider the case that $\varphi$ is bicyclic with elementary divisors $(x-1)^{2t}$ and nonhyperbolic. By \ref{SKEWINV8}, $\varphi$ has a
	symplectic square root $\psi$ with minimal polynomial $(x^2+1)^n$.
	By \ref{SKEWINV2} and \ref{SKEWINV3}, $\psi$ is reversed by a symplectic skew-involution $\eta$. Hence $\varphi^{\eta} = (\psi^2)^{\eta} = (\psi^{\eta})^2 = (\psi^{-1})^2 = \varphi^{-1}$.
\end{proof}

\section{Products of an involution and a skew-involution}
\mbox{}

\begin{lemma}                                                           \label{INV_SKEWINV1}
	Let $P \in \GL(n, K)$. Assume that $P^2$ has no elementary divisors $(x+1)^t$ of odd degree. Then $P$ is similar to $-P^{-1}$ if and only if  $P$ is a product of an involution and a skew-involution.
\end{lemma}

\begin{proof}
	If $P = SH$, where $S^2 = \Idm_n = - H^2$, then $P^S = HS = -H^{-1}S = -(SH)^{-1}$. Conversely, by \ref{SKEWINV0}, $P$ is similar to a matrix
	\[
	Q = \left (\begin{array} {cc} 0 & \Idm_m  \\  \Idm_m & D \end{array} \right )
	= \left (\begin{array} {cc} \Idm_m & 0  \\ D & - \Idm_m\end{array} \right )
	\left (\begin{array} {cc} 0 & \Idm_m  \\  -\Idm_m & 0 \end{array} \right ).
	\]	
\end{proof}

\begin{corollary}                                             \label{INV_SKEWINV2}
	Let $P \in \SpG(2n, K)$. If 
	\[
	P = \left (\begin{array} {cc} A & 0\\ 0 & A^+ \end{array} \right ),
	\]
	where $A$ similar to $-A^{-1}$ or $-A$. 
\end {corollary}

\begin{proof}
	 If $A$ is similar to $-A^{-1}$, apply \ref{INV_SKEWINV1}.
	 So let $A$ be similar to $-A$. Then we may assume that
	 \[
	 A = 
	 \left (\begin{array} {cc} 0 & \Idm_m\\ C & 0 \end{array} \right ).
	 \]
	 Let $C = RT$ is the product of 2 symmetric matrices $R$ and $T$.
	 Then 
	 \[
	 A = \left (\begin{array} {cc} 0 & -R\\ R & 0 \end{array} \right )
	 \left (\begin{array} {cc} T & 0\\ 0 & -T^{-1} \end{array} \right )
	 \]
	 is the product of an anti-symmetric matrix $H$ and a symmetric matrix $S$ and 

	 \[
	 P = \left (\begin{array} {cc} A & 0\\ 0 & A^+ \end{array} \right )
	 = \left (\begin{array} {cc} 0 & H\\ -H^+ & 0 \end{array} \right )
	 \left (\begin{array} {cc} 0 & -S^{-1}\\ S & 0 \end{array} \right ).
	 \]
	 is a product of a symplectic involution and a skew-involution.
\end{proof}

\begin{lemma}                                             \label{INV_SKEWINV3}
	Let $\varphi \in \SpG(V)$ be cyclic with minimal polynomial 
	\begin{enumerate}
	   \item $p(x)^t$, where $p(x)$ is irreducible and $p(x) = p(-x)(-1)^{\partial p}$, or
	   \item $p(x)^t q(x)^t$, where $p(x)$ is irreducible and prime to $q(x) =  p(-x)(-1)^{\partial p}$.
	\end{enumerate}
	Then there exists an involution 
	$\sigma \in \SpG(V)$ such that $\varphi^{\sigma} = -\varphi^{-1}$ if and only if $\varphi^2$ is hyperbolic.
\end{lemma}
	
\begin{proof}
	Clealy, if $\varphi^{\sigma} = -\varphi^{-1}$, then $\varphi^2$ is bireflectional, and by \ref{theorem-2}, $\varphi^2$ is hyperbolic.
		Conversely, let $U = \langle u \rangle_{\varphi^2}$ be a totally degenerate $\varphi^2$-cyclic subspace of $V$ with dimension $\frac{1}{2} \dim V$. Then $V = \langle u \rangle_{\varphi}$: 
		
		This is clear if $\varphi$ is orthogonally indecomposable. 
		So let $\mu(\varphi) = p(x)^t q(x)^t$, where $p(x)$ is prime to $q(x)$. The only maximal $\varphi$-invariant subspaces are 
		$\ker p(\varphi)^{t-1} \oplus \ker p(-\varphi)^t$ and 
		$\ker p(-\varphi)^{t-1} \oplus \ker p(\varphi)^t$.
		Replacing $\varphi$ with $-\varphi$, we may assume  that 
		$u \in \ker p(\varphi)^{t-1} \oplus \ker p(-\varphi)^t$.
		Then $\ker p(\varphi) \le U$ and $u p(\varphi)^{-1} p(-\varphi)^{t-1} \in \ker p(-\varphi)$. Hence 
		$\ker p(\varphi) \oplus \ker p(-\varphi) \le U$, a contradiction.
		
		Define $\sigma$ by
		$u\varphi^i \sigma = (-1)^i u \varphi^{-i}$. Then $\sigma^2 = 1$ and
		$\varphi^{\sigma} = -\varphi^{-1}$. Finally, $\sigma$ is symplectic, since $U$ is totally degenerate.
	\end{proof}

 \begin{lemma}                                 \label{INV_SKEWINV4}
	Let $\varphi \in \SpG(V)$ be cyclic with
	minimal polynomial $(x^2+1)^n$, where $n$ is odd.
	Then  $\varphi$ is the product of a symplectic involution and a skew-involution.
\end{lemma}
	
\begin{proof}
		$\varphi^2$ must be orthogonally indecomposable of type 1. By \ref{prop-1}, $\varphi^2$ is hyperbolic. Apply \ref{INV_SKEWINV3}.
\end{proof}

\begin{lemma}  \label{INV_SKEWINV5}                       
	Let $K = \GF(q)$ be a finite field with $q \equiv 3 \mod 4$, 
	and let $\varphi \in \SpG(V)$. Assume 
	that $\Fix(\varphi^2) = 0$.
	Then the following are equivalent
	\begin{enumerate}
		\item $\varphi^{\sigma} = -\varphi^{-1}$ for some involution $\sigma \in \SpG(V)$;
		\item $\varphi$ is similar to $-\varphi$ and every elementary divisor $e(x) \in \{(x^2+1)^{2t}; t \in \mathbb{N}\}$
		of $\varphi$ occurs with even multiplicity.
	\end{enumerate}
\end{lemma}

\begin{proof}
	$1 \rightarrow 2$: Clearly, $(\varphi^2)^{\sigma} = (\varphi^2)^{-1}$. Let  $(x^2+1)^{2t}$
	be an elementary divisor of $\varphi$ with multiplicity $s$. 
	Let $W = \ker (\varphi^2+1)^{t+1} \cap V(\varphi^2+1)^{2n - t-1}$. Then 
	$\varphi$ induces a homocyclic isometry with minimal polynomial
	$(x^2+1)^2$ on $W/\rad W$.
	We have $\dim W/\rad W = 4s$, so by induction, we may assume $V = W$. Clearly, $\sigma$ inverts $\varphi^2$.
	It follows from \ref{SYMPINV0} that $\dim \Fix(\sigma) = \frac{\dim V}{2} = 2s$. Now $\ker (\varphi^2+1)$ is $\sigma$-invariant and totally isotropic. Hence $\dim \Fix(\sigma) \cap \ker (\varphi^2+1) = s$. 
	So $s$ must be even, as $\sigma$ is the product of 2 skew-involutions on $\ker (\varphi^2+1)$.
	
	$2 \rightarrow 1$:
	We may assume that either
	\begin{itemize}
		\item $\varphi$ is cyclic with minimal polynomial $(x^2+1)^{2m+1}$, or
		\item $\varphi$ is cyclic with minimal polynomial prime to $x^2+1$, or
		\item $\varphi$ is bicyclic with elementary divisors $e_1 = e_2 = (x^2+1)^{2m}$.
	\end{itemize}
	Let $\varphi$ be  cyclic. If $\mu_{\varphi} = (x^2+1)^{2m+1}$ we can apply \ref{INV_SKEWINV4}.
	
	Let $\mu_{\varphi} = p(x)^t$, where $p$ is prime  and $p(x)  \ne x^2+1$. 
	There exists a polynomial $q \in K[x]$ such that $p(x) = q(x^2)$.
	Hence $\varphi^2$ is bicyclic with elementary divisors $e_1 = e_2 = q(x)^t$.
	Since $q(x) \ne x \pm 1$, it follows that $\varphi^2$ is hyperbolic by \ref{WALL1}. Apply \ref{INV_SKEWINV3}.
	
	If $\varphi$ is bicyclic with elementary divisors $e_1 = e_2 = (x^2+1)^{2m}$, then $\varphi$ is hyperbolic by \ref{WALL1},
	and we can apply \ref{INV_SKEWINV2}.
\end{proof}

\begin{lemma}                                                                                \label{INV_SKEWINV6}                     
	Let $K = \GF(q)$ be a finite field with $q \equiv 3 \mod 4$, 
	and let $\varphi \in \SpG(V, f)$. Suppose that $\varphi^2$ is unipotent.
	Then the following are equivalent:
	\begin{enumerate}
		\item $\varphi^{\sigma} = -\varphi^{-1}$ for some involution $\sigma \in \SpG(V)$;
		\item $\varphi$ is conjugate to $-\varphi^{-1}$;
		\item $\varphi$ has a decomposition $\varphi = \varphi_1 \perp -\varphi_1 \perp \dots \perp \varphi_s \perp  -\varphi_s$, 
		where $\varphi_j$ is orthogonally indecomposable of type 1 or 2
		and $(\varphi_j \perp -\varphi_j)^2$ is hyperbolic.
	\end{enumerate}
\end{lemma}

\begin{proof}
	$1 \rightarrow 2$: clear\\
	$2 \rightarrow 3$: 
	Let $\alpha \in \SpG(V)$ with $\varphi^{\alpha} = -\varphi^{-1}$.
	Then $V$ has a decomposition $V = U_1 \perp U_1\alpha \perp \dots \perp U_s \perp U_s\alpha$, where the subspaces 
	$U_1, \dots, U_s$ are $\varphi$-invariant and orthogonally indecomposable of type 1 or 2.
	
	If $U_j$ is of type 1, then already  $\varphi_j$ is hyperbolic by \ref{prop-1}. So assume that $U:=U_j$ is of type 2.
	Let $U = \langle u \rangle_{\varphi}$. Then $\langle u + u\alpha\rangle_{\varphi^2}$ and 
	$\langle u - u\alpha\rangle_{\varphi^2}$ are totally isotropic:
	We have $f(u \alpha \varphi^{2t}, u\alpha) = f(u \alpha, u \alpha \varphi^{-2t}) =  f(u \alpha, u \varphi^{2t} \alpha) = f(u, u \varphi^{2t})$.
	
	Hence 
	$f((u \pm u\alpha)\varphi^{2t}, u \pm u\alpha) = f(u\varphi^{2t}, u) + f(u\alpha\varphi^{2t}, u\alpha) = 0$.
	Furthermore, $U \perp U\alpha = \langle u + u\alpha\rangle_{\varphi^2} \oplus \langle u - u\alpha\rangle_{\varphi^2}$.
	
	$3 \rightarrow 1$: We may assume that $s=1$. Apply \ref{INV_SKEWINV2} and \ref{INV_SKEWINV3}.
\end{proof}

\begin{proof}[Proof of theorem \ref{theorem-5}]
	Apply \ref{INV_SKEWINV5} and  \ref{INV_SKEWINV6}.
\end{proof}


\ifdraft{\listoflabels}{}

\end{document}
\typeout{get arXiv to do 4 passes: Label(s) may have changed. Rerun}